\documentclass[12pt]{article}
\usepackage{amsmath}
\usepackage{amssymb}
\usepackage{graphicx}
\usepackage[cp1251]{inputenc}
\usepackage[russian]{babel}
\newcommand{\il}[2]{\int\limits_{#1}^{#2}}

\newcommand{\ph}{\phantom{a}}
\newcommand{\phh}{\phantom{aaa}}

\newcommand{\sist}[2]{\left\{
\begin{array}{l}
{#1}\\
\ph\\
{#2}
\end{array}
\right.}

\begin{document}

\vskip 20pt

MSC 34C10

\vskip 20pt

\centerline{\bf Oscillation and non oscillation criteria for linear  nonhomogeneous}
\centerline{\bf   systems of two first-order ordinary differential equations}

\vskip 20 pt

\centerline{\bf G. A. Grigorian}
\centerline{\it Institute  of Mathematics NAS of Armenia}
\centerline{\it E -mail: mathphys2@instmath.sci.am}
\vskip 20 pt

\noindent
Abstract. The Riccati equation method is used to establish an oscillatory and a non oscillatory criteria for nonhomogeneous linear systems of two first-order ordinary \linebreak differential equations. It is shown that the obtained oscillatory criterion is a generalization of J. S. W. Wong's oscillatory criterion.

\vskip 20 pt


\vskip 20 pt

{\bf 1. Introduction.} Let $p(t), \ph q(t), \ph r(t), \ph s(t), \ph f(t)$ and $g(t)$ be real-valued continuous functions on $[t_0,+\infty)$. Consider the linear nonhomogeneous systems of differential \linebreak   equations
$$
\sist{\phi' = p(t)\phi + q(t) \psi + f(t),}{\psi' = r(t)\phi + s(t)\psi + g(t), \phh t \ge t_0.} \eqno (1.1)
$$
and the corresponding homogeneous one
$$
\sist{\phi' = p(t)\phi + q(t) \psi,}{\psi' = r(t)\phi + s(t)\psi, \phh t \ge t_0.} \eqno (1.2)
$$

{\bf Definition 1.1.} {\it The system (1.1) ((1.2)) is called oscillatory if for each of its  solutions $(\phi(t), \psi(t))$ the function  $\phi(t)$ has arbitrary large zeroes, otherwise it is called non \linebreak oscillatory.}

{\bf Definition 1.2.} {\it  The system (1.1) ((1.2)) is called oscillatory on the interval $[t_1,t_2] \linebreak (\subset [t_0,+\infty))$  if for each of  its  solutions $(\phi(t), \psi(t))$ the function  $\phi(t)$ has a zero on $[t_1,t_2]$.}

Let $a(t), \ph b(t), \ph c(t)$ and $d(t)$ be real-valued continuous function on $[t_0,+\infty)$, and let $a(t) > 0, \ph t \ge t_0$. Along with the systems (1.1) and (1.2) consider the second order linear non homogeneous equation
$$
(a(t)\phi')' + b(t) \phi' + c(t) \phi = d(t), \phh t \ge t_0. \eqno (1.3)
$$
The substitution
$$
\psi = a(t) \phi', \phh t \ge t_0. \eqno (1.4)
$$
reduces Eq. (1.3) to the following linear system
$$
\sist{\phi' = \frac{1}{a(t)} \psi,} {\psi' = - c(t) \phi - \frac{b(t)}{a(t)} \psi + d(t), \ph t \ge t_0.} \eqno (1.5)
$$

{\bf Definition 1.3.} {\it Eq. (1.3) is called oscillatory if eacch of its solutions has arbitrary large zeroes, otherwise it is called non oscillatory.}

{\bf Definition 1.4.} {\it Eq. (1.3) is called oscillatory  on the interval $[t_1,t_2] (\subset [t_0,+\infty))$ if each of its  solutions has a zero on  $[t_1,t_2]$.}

Obviously Eq. (1.3) is oscillatory (non oscillatory, oscillatory on the interval $[t_1,t_2] \linebreak (\subset [t_0,+\infty))$) if and only if the system (1.5) is  oscillatory (non oscillatory, oscillatory on the interval $[t_1,t_2] (\subset [t_0,+\infty))$).

Note that the system (1.5) is a particular case of the system (1.1) therefore the system (1.1) can be considered as a generalization of Eq. (1.3). The system (1.1) (especially  Eq. (1.3)) has an important theoretical interest and practical applications. Therefore study the question  of oscillation and non oscillation of the system (1.1) (in particular of Eq. (1.3)) is an important problem of qualitative theory of differential equations. To study of this question for Eq. (1.3) are devoted many works (see [1 - 10] and cited works therein). Among them notice the El-Saedy's oscillation theorem for undamped ($b(t)\equiv 0$) Eq. (1.3) (see [1], or [2, Theorem A]). Combining the Riccati equation method with  a variational technique  J. S. W. Wong obtained the following generalization of El-Saedy's result (see [2, Theorem 1]).

{\bf Theorem 1.1.} {\it Suppose that $b(t)\equiv 0$ and for any $T \ge t_0$ there exist $T \le s_1 < t_1 \le s_2 < t_2$ such that
$$
d(t)\sist{\le 0, \ph t\in [s_1, t_1],}{\ge 0, \ph t\in [s_2,t_2].}
$$
Denote $D(s_i,t_i) \equiv \{u \in \mathbb{C}^1[s_i,t_i]: u(t)\not\equiv 0, \ph u(s_i) = u(t_i) = 0\}, \ph i =1,2.$. If there exists $u \in D(s_i,t_i)$ such that
$$
\il{s_i}{t_i}\bigl(c(\tau) u^2(\tau) - a(\tau) u'(\tau)^2\bigr) d\tau \ge 0 \eqno (1.6)
$$
for $i=1,2,$ then Eq. (1.3) is oscillatory.}

Later using the same approach of the proof  of Theorem 1.1  Q. Kong and M. Pasic extended the J. S. W. Wong's result to the damped Eq. (1.3) (see [3, Theorem 2]).
Unfortunately despite  the significance of these results their conditions contain free para-\linebreak meter-functions. This makes it hard to further application them to the concrete equations (despite the fact that Q. Kong and M. Pasic found in  [3] an one parameter family of equations satisfying the conditions of Theorem 1.1).

Unlike of the approach used for proving mentioned above oscillatory criteria of works [2,3]  in this paper we use only the Riccati equation method for proving an oscillatory and a non oscillatory criteria for the system (1.1). We show that the obtained oscillatory criterion is a generalization of Theorem 1.1.

{\bf 2. Auxiliary  propositions.} Let $f_k(t), \ph g_k(t), \ph h_k(t), \ph k=1,2$ be real-valued conti-\linebreak nuous functions on $[t_0,+\infty)$. Consider the Riccati equations
$$
y' + f_k(t) y^2 + g_k(t) y + h_k(t) = 0, \phh t\ge t_0. \eqno (2.1_k)
$$
$k=1,2$ and the differential inequalities
$$
\eta + f_k(t) \eta^2 + g_k(t) \eta + h_k(t) \ge 0, \phh t\ge t_0. \eqno (2.2_k)
$$
$k=1,2$.

{\bf Remark 2.1.} {\it Every solution  of Eq. $(2.1_2)$ on $[t_0,t_1)$ is also a solution of the inequality $(2.2_2)$ on $[t_0,t_1)$.}

{\bf Remark 2.2.} {\it If $f_1(t) \ge 0, \ph t\in [t_0,t_1)$, then every solution of the linear equation
$$
\zeta' + g_1(t)\zeta + h_1(t) = 0, \phh t\in [t_0,t_1)
$$
is also a solution of the inequality $(2.2_1)$ on $[t_0,t_1)$.}

{\bf Theorem 2.1.} {\it  Let $y_2(t)$ be a solution of Eq. $(2.1_2)$ on $[t_0,\tau_0) \ph (t_0 < \tau_0 \le +\infty)$ and let $\eta_1(t)$ and $\eta_2(t)$ be solutions of the inequalities $(2.2_1)$ and $(2.2_2)$  respectively on $[t_0,\tau_0)$ such that $y_2(t_0) \le \eta_k(t_0) \ph k=1,2.$ In addition let the following conditions be satisfied: $f_1(t) \ge 0, \ph \gamma - y_2(t_0) + \il{t_0}{t}\exp\biggl\{\il{t_0}{\tau}[f_1(s)(\eta_1(s) + \eta_2(s)) + g_1(s)]ds\biggr\}\biggl[(f_2(\tau) - f_1(\tau))^2 y_2^2(\tau) + (g_2(\tau) - g_1(\tau)) y_2(\tau) + h_2(\tau) - h_1(\tau)\biggr] d \tau \ge 0, \ph t\in [t_0,\tau_0)$ for some $\gamma \in [y_2(t_0), \eta_1(t_0)]$. Then Eq. $(2.1_1)$ has a solution $y_1(t)$ on $[t_0,\tau_0)$ with $y_1(t_0) \ge \gamma$ and $y_1(t) \ge y_2(t), \ph t\in [t_0,\tau_0)$.}

See the proof in [11].

For any $\lambda \in \mathbb{R}$ set:
$$
\alpha_\lambda(t) \equiv \exp\biggl\{\il{t_0}{t}p(\tau)d\tau\biggr\}\biggl[\lambda + \il{t_0}{t} \exp\biggl\{-\il{t_0}{\tau}p(s)d s\biggr\}f(\tau)d\tau\biggr], \phh t\ge  t_0
$$
and in the system (1.1) substitute
$$
\phi = \phi_1 + \alpha_\lambda(t). \eqno (2.3)
$$
Since $\alpha_\lambda(t)$ is a solution of the equation
$$
\alpha' = p(t)\alpha + f(t), \phh t \ge t_0
$$
we obtain
$$
\sist{\phi_1' = p(t)\phi_1 + q(t) \psi,}{\psi' = r(t)\phi_1 + s(t) \psi + g_\lambda(t),\ph t \ge t_0,} \eqno (2.4)
$$
where $g_\lambda(t)\equiv r(t)\alpha_\lambda(t) + g(t), \ph t \ge t_0.$ In the obtained system substitute
$$
\psi = y \phi_1.
$$
We get
$$
\sist{\phi_1' = [p(t) + q(t) y] \phi_1,}{y'\phi_1  + [q(t) y^2 + E(t) y - r(t)]\phi_1 = g_\lambda(t), \ph t \ge t_0,}
$$
where $E(t)\equiv p(t) - s(t), \ph t\ge t_0.$ It follows from here that if $(\phi_1(t), \psi(t))$ is a solution of the system (2.4) such that $\phi_1(t) \ne 0, \ph t\in [t_1,t_2) \ph (t_0\le t_1 < t_2 \le +\infty)$ then
$$
\phi_1(t) = \phi_1(t_1)\exp\biggl\{\il{t_0}{t}\Bigl[p(\tau) + q(\tau) y(\tau)\Bigr]d \tau\biggr\}, \ph \psi(t) = y(t)\phi_1(t), \ph t\in [t_1,t_2), \eqno (2.5)
$$
where $y(t)$ is the solution of the Riccati equation
$$
y' + q(t)y^2 + E(t) y - h_{\phi_1}(t) = 0, \phh t \in [t_1,t_2),
$$
with $y(t_1) = \frac{\psi(t_1)}{\phi_1(t_1)},$ where $h_{\phi_1}(t) \equiv r(t)+ \frac{g_\lambda(t)}{\phi_1(t)}, \phh t\in [t_1,t_2)$.

{\bf 3. Main results.} In this section we prove a non oscillatory and an oscillatory criteria for the system (1.1). We show that the obtained oscillatory criterion is a generalization of Theorem 1.1.

{\bf Theorem 3.1.} {\it Let the following conditions be satisfied:

\noindent
1) $q(t) \ge 0, \phh t \ge t_0;$

\noindent
2) there exists $\lambda \ge 0$ such that

\noindent
2$_1$) $\alpha_\lambda(t) \ge 0, \phh t \ge t_0$;

\noindent
2$_2$) $ r(t)\alpha_\lambda(t) + g(t) \ge 0, \phh t \ge t_0.$

\noindent
Then if the system (1.2) is non oscillatory the system (1.1) is also non oscillatory.}

Proof. Assume the system (1.2) is non oscillatory. Then there exists its  a solution $(\phi_0(t), \psi_0(t))$  such that $\phi_0(t) \ne 0, \ph t \ge t_1$ for some $t_1 \ge t_0$. It follows from here that $y_0(t) \equiv \frac{\psi_0(t)}{\phi_0(t)}, \ph t \ge t_1$ is a solution of the Riccati equation
$$
y' + q(t) y^2 + E(t) y - r(t) = 0, \phh t \ge t_0. \eqno (3.1)
$$
on $[t_1, +\infty)$. Let $(\phi(t), \psi(t))$ be a solution of the system (1.1) with $\phi(t_1) > \alpha_\lambda(t_1)$,
$$
\psi(t_1) = (\phi(t_1) - \alpha_\lambda(t_1)) y_0(t_1) \eqno (3.2)
$$
and let $\phi_1(t) \equiv \phi(t) - \alpha_\lambda(t), \ph t \ge t_0$. Then
$$
\phi_1(t_1) > 0 \eqno (3.3)
$$
and by (2.3) $(\phi_1(t), \psi(t))$ is a solution of the system (2.4). Show that
$$
\phi_1(t) > 0, \phh t \ge t_1. \eqno (3.4)
$$
Suppose it is not so. Then by (3.3) there exists $T > t_1$ such that $\phi_1(t) > 0, \ph t \in [t_1,T)$ and
$$
\phi_1(T) = 0. \eqno (3.5)
$$
In virtue of the first equality of (2.5) we have
$$
\phi_1(t) = \phi_1(t_1)\exp\biggl\{\il{t_1}{t}\Bigl[p(\tau) + q(\tau) y(\tau)\Bigr]d \tau\biggr\}, \phh t \in [t_1, T), \eqno (3.6)
$$
where $y(t)$ is a solution of the Riccati equation
$$
y' + q(t) y^2 + E(t) y - r(t) - \frac{r(t)\alpha_\lambda(t) + g(t)}{\phi_1(t)} = 0, \phh t \in [t_1,T). \eqno (3.7)
$$
By the second equality of (2.5) from (3.2) it follows that $y(t_1) = y_0(t_1).$ Then applying Theorem 2.1 to the pair of Equations (3.1) and (3.7) on the basis of the conditions 1) and 2$_2$) we conclude that $y(t) \ge y_0(t), \ph t \in [t_1,T)$. This together with 1) and (3.6) implies that
$$
\phi_1(T) \ge \phi_1(t_1)\exp\biggl\{\il{t_1}{T}[p(\tau) + q(\tau) y(\tau)]d \tau\biggr\} > 0,
$$
which contradicts (3.5). The obtained contradiction proves (3.4). It follows from the condition 2$_1$) that
$$
\phi(t) = \phi_1(t) + \alpha_\lambda(t) \ge \phi_1(t), \phh t \ge t_1.
$$
This together with (3.4) implies that $\phi(t) > 0, \ph t \ge t_1$. Therefore the system (1.1) is non oscillatory. The theorem is proved.

{\bf Remark 3.1.} {\it Some explicit non oscillatory criteria for the system (1.2) are obtained in [12 - 14].}

{\bf Theorem 3.2.} {\it Let the following conditions be satisfied:

\noindent
1) $q(t) \ge 0, \ph t \ge t_0;$

\noindent
3) for every $T > t_0$ there exist $T \le s_1 <t_1\le s_2 < t_2$ and $\lambda \in \mathbb{R}$ such that

\noindent
3$_1$) $(-1)^k\alpha_\lambda(t) \ge 0, \ph t\in[s_k,t_k], \ph k=1,2;$

\noindent
3$_2$)  $(-1)^k[r(t)\alpha_\lambda(t) + g(t)] \ge 0, \ph t\in[s_k,t_k], \ph k=1,2;$

\noindent
3$_3$) the system (1.1) is oscillatory on the intervals $[s_k, t_k], \ph k=1,2.$.

\noindent
Then the system (1.1) is oscillatory.}

Proof. Suppose the system (1.1) is not oscillatory. Then there exists its a solution $(\phi(t), \psi(t))$ such that $\phi(t) \ne 0, \ph t \ge T,$ for some $T \ge t_0$. We set: $k=1$ if $\phi(t) > 0, \ph t \ge T$ and $k=2$ if $\phi(t) < 0, \ph t \ge T$. Let $\phi_1(t)\equiv \phi(t) - \alpha_\lambda(t), \ph t \ge t_0.$ Then by (2.3) $(\phi_1(t), \psi(t))$ is a solution of the system (2.4) and by 3$_1$) we have $\phi_1(t) \ne 0, \ph t\in [s_k, t_k +\varepsilon]$ for some $\varepsilon > 0$, so by virtue of (2.5)
$$
\phi_1(t) = \phi_1(s_1)\exp\biggl\{\il{s_k}{t}\Bigl[p(\tau) + q(\tau) y(\tau)\Bigr]d \tau\biggr\},\phh t\in [s_k,t_k],
$$
where $y(t)$ is a solution of the Riccati equation
$$
y' + q(t) y^2 + E(t) y -  r(t) - \frac{r(t)\alpha_\lambda(t) + g(t)}{\phi_1(t)} = 0, \phh t \in [s_k, t_k+\varepsilon) \eqno (3.8)
$$
with $y(s_k) = \frac{\psi(s_k)}{\phi_1(s_k)}.$ Consider the Riccati equation
$$
y' + q(t) y^2 + E(t) y -  r(t) = 0, \phh t \in [s_k, t_k+\varepsilon) \eqno (3.9)
$$
From 3$_2$) it follows that
$$
\frac{r(t)\alpha_\lambda(t) + g(t)}{\phi_1(t)} \le 0, \phh t\in [s_k, t_k].
$$
Then applying Theorem 2.1 to the pair of equations (3.8) and (3.9) we conclude that Eq. (3.9) has a solution $y_1(t)$ on $[s_k,t_k)$ with $y_1(s_1) \ge y(s_1)$ and    $y_1(t) \ge y(t), \ph t\in [s_k,t_k)$. Hence according to (2.5) the pair of functions
$$
\phi_2(t) \equiv \exp\biggl\{\il{t_0}{t}[p(\tau) + q(\tau) y(\tau)]d \tau\biggr\}, \ph \psi_2(t) \equiv y_1(t)\phi_2(t), \ph t\in [s_k,t_k]
$$
form a solution $(\phi_2(t), \psi_2(t))$ of the system (1.2) on $[s_k, t_k)$ and is continuable on $[s_k,t_k]$ as a solution of the system (1.2), for which $\phi_2(t) \ne 0, \ph t \in [s_k, t_k]$. Therefore the system (1.2) is not oscillatory on $[s_k, t_k]$, which contradicts the condition 3$_3$). The obtained contradiction completes the proof of the theorem.

{\bf Remark 3.2.} {\it An explicit interval oscillatory criterion for the system (1.2)  is obtained in [15] (see also [16]). An interval oscillatory criteria for unforced Eq. (1.3) is obtained in [17].}

One can easily show (using the same way of proof of Theorem 1.1) that under the conditions of Theorem 1.1 the equation
$$
(a(t)\phi')' + c(t)\phi = 0, \phh t \ge t_0  \eqno (3.10)
$$
is oscillatory on the intervals $[s_k,t_k], \ph k=1,2$ for $s_k, t_k, \ph k=1,2$ indicated in the conditions of Theorem 1.1. on the other hand in virtue of the connection (1.4) between Eq. (1.3) and the  system (1.5) the conditions of Theorem 3.2 for the system (1.5) in the particular case $\lambda = 0$ can be reduced (equivalently translated) for undamped Eq. (1.3) to the following ones.

\noindent
For every $T\ge t_0$ there exist $T\le s_1 < t_1 \le s_2 < t_2$ such that
$$
d(t)\sist{\le 0, \ph t\in [s_1,t_1],}{\ge 0, \ph t\in [s_2,t_2];}
$$
the equation (3.10) is oscillatory on the intervals $[s_k, t_k], \ph k=1,2$.

\noindent
Therefore Theorem 3.2 is a generalization of Theorem 1.1, moreover the condition 3$_3$) of Theorem 3.2 is preferable to the condition (1.6), since each explicit oscillatory criterion for the system (1.2) generates according to Theorem 3.2 an explicit oscillatory criterion, whereas (1.6) does not have this property (since no explicit conditions (or an explicit condition) on the coefficients of Eq. (3.10) that ensure (1.6) have been specified yet).

\pagebreak

\centerline{\bf References}

\vskip 20pt

\noindent
1. M. A. El-Sayed, An oscillation criterion for a forced second order linear equations, Proc.  \linebreak \phantom{a}  Amer. Math. Soc. 118 (19930, 813 - 817.

\noindent
2. J. S. W. Wong, Oscillation Criteria for a Forced Second-Order Linear Differential \linebreak \phantom{a}  Equations, J. Math. Anal. Appl., 231, 235 - 240 (1999).

\noindent
3. Q. Kong and M. Pasic,   Second - Order Differential Equations: Some Significant results \linebreak \phantom{a}   Due to James S. W. Wong,   Diff Eq. and Appl., vol. 6, Num. 1 (2.14), 99 - 163.

\noindent
4. M. K. Kwong, J. S. W. Wong, On the oscillation of Hill's equation under periodic  \linebreak \phantom{a} forcing, J. Math.   Anal. Appl., 320, (2006) 37 - 55.

\noindent
5. A. Skidmore and W. Leighton, On the differential equation $y'' + p(x) y = f(x),$, J. Math. \linebreak \phantom{a} Anal. Appl., 43 (1973) 46 -55.

\noindent
6. A. Skidmore and J.J. Bowers, Oscillation behavior of $y'' + p(x) y = f(x),$, J. Math. \linebreak \phantom{a} Anal. Appl., 49 (1975), 317 - 323.

\noindent
7. S. M. Rainkin, Oscillation theorems for second order nonhomogeneous linear differential \linebreak \phantom{a} equations, J. Math. Anal. Appl., 53 (1976) 550 - 553.

\noindent
8. A. G. Kortsatas, Maintenance of oscillation under the effect of a periodic forcing term,\linebreak \phantom{a} Proc. Amer. Math. Soc. 33 (1972), 377 - 383.

\noindent
9. M. S. Kenner, Solutions of a certain linear nonhomogeneous second order differential \linebreak \phantom{a} equations, Appl. Anal. 1. (1971) 57 - 63.

\noindent
10. Y. G. Sun, C. H. Ou, and J. S. W. Wong, Interval Oscillation Theorems for a \linebreak \phantom{a} Second-Order Linear differential equations, Compurters and mathematics with \linebreak \phantom{a} Applications, 48, (2004) 1693 - 1699.

\noindent
11. G. A. Grigorian,  On two comparison tests for second-order linear  ordinary\linebreak \phantom{aa} differential equations (Russian) Differ. Uravn. 47 (2011), no. 9, 1225 - 1240; trans-\linebreak \phantom{aa} lation in Differ. Equ. 47 (2011), no. 9 1237 - 1252, 34C10.

\noindent
12. G. A. Grigorian. Oscillatory Criteria for the Systems of Two First - Order Linear\linebreak \phantom{aa} Ordinary Differential equations. Rocky Mountain Journal of Mathematics, vol. 47,\linebreak \phantom{aa} no. 5, 2017, pp. 1497 - 1524.

\noindent
13	G. A. Grigoryan, Two comparison criteria for scalar Riccati equations and their \linebreak \phantom{a}  applications, Izv. Vyssh. Uchebn. Zaved. Mat., 2012, 11,  20–35; Russian Math.\linebreak \phantom{aa} (Iz. VUZ), 56:11 (2012), 17–30.

\noindent
14.	G. A. Grigoryan, Criteria of global solvability for Riccati scalar equations, Izv. Vyssh. \linebreak \phantom{a} Uchebn. Zaved. Mat., 2015, 3,  35–48; Russian Math. (Iz. VUZ), 59:3 (2015), 31–42.

\noindent
15. G. A. Grigorian, Oscillatory and non oscillatory criteria for the systems of two linear \linebreak \phantom{a} first order two by two dimensional matrix ordinary differential equations. Arch. Math.\linebreak \phantom{a}  54 (2018), no. 4, 189–203.

\noindent
16. G. A. Grigorian.   Interval oscillation criteria for linear matrix Hamiltonian systems,\linebreak \phantom{a} vol. 50 (2020), No. 6, 2047–2057.

\noindent
17.	G. A. Grigoryan, Some properties of solutions of second-order linear ordinary differential \linebreak \phantom{aa} equations, Trudy Inst. Mat. i Mekh. UrO RAN, 19:1 (2013),  69–80.

\end{document}